\documentclass[11pt, draft]{amsart}
\usepackage{amssymb}
\usepackage[mathscr]{eucal}

\theoremstyle{plain}
\newtheorem*{mainthm}{Main Theorem}
\newtheorem{theorem}{Theorem}

\newtheorem*{lemma}{Lemma}
\newtheorem{corollary}[theorem]{Corollary}

\theoremstyle{remark}
\newtheorem*{remark}{Remark}

\newcommand\C{\mathbb C}
\newcommand\f{{\underline{f}}}
\renewcommand\H{{\underline{H}}}
\renewcommand\ker{\operatorname{ker}}
\newcommand\K{\mathbb K}
\newcommand\la{\langle}
\newcommand\rng{\operatorname{rng}}
\newcommand\ot{\otimes}
\newcommand\R{\mathbb R}
\newcommand\ra{\rangle}
\newcommand\tphi{\tilde \phi}
\newcommand\tr{\operatorname{tr}}
\newcommand\x{{\underline{x}}}
\newcommand\X{{\underline{X}}}
\newcommand\Xs{{\underline{X}'}}
\newcommand\y{{\underline{y}}}

\begin{document}
\title[]{Orthogonality preserving transformations on indefinite inner
product spaces:\\
Generalization of Uhlhorn's version of Wigner's theorem}
\author{LAJOS MOLN\'AR}
\address{Institute of Mathematics and Informatics\\
         University of Debrecen\\
         4010 Debrecen, P.O.Box 12, Hungary}
\email{molnarl@math.klte.hu \newline \phantom{M\,}
{\normalfont {\it URL}:\,\, http://neumann.math.klte.hu/\~{}molnarl/}}
\thanks{  This research was supported by the
          Hungarian National Foundation for Scientific Research
          (OTKA), Grant No. T030082, T031995, and by
          the Ministry of Education, Hungary, Reg.
          No. FKFP 0349/2000}
\date{\today}
\begin{abstract}
We present an analogue of Uhlhorn's version of Wigner's theorem on
symmetry transformations for the case of indefinite inner product
spaces. This significantly generalizes a result of Van den Broek.
The proof is based on our main theorem, which describes the form
of all bijective transformations on the set of all rank-one
idempotents of a Banach space which preserve zero products in both
directions.
\end{abstract}
\maketitle

\section{Introduction and statement of the results}

Wigner's theorem on symmetry transformations plays a fundamental
role in quantum mechanics. It states that any quantum mechanical
invariance transformation (symmetry transformation) can be
represented by a unitary or antiunitary operator on a complex
Hilbert space and that, conversely, any operator of that kind
represents an invariance transformation. In mathematical language,
the result can be reformulated in the following way. If $H$ is a
complex Hilbert space and $T$ is a bijective transformation on the
set of all 1-dimensional linear subspaces of $H$ which preserves
the angle between every pair of such subspaces (in the terminology
of quantum mechanics, this angle is called a transition
probability), then $T$ is induced by either a unitary or an
antiunitary operator $U$ on $H$. This means that for every
1-dimensional subspace $L$ of $H$ we have $T(L)=U[L]=\{ Ux \, :\,
x\in L\}$. In his famous paper \cite{Uhlhorn}, Uhlhorn generalized
this result by requiring only that $T$ preserves the orthogonality
between the 1-dimensional subspaces of $H$. This is a significant
achievement since Uhlhorn's transformation preserves only the
logical structure of the quantum mechanical system in question
while Wigner's transformation preserves its complete probabilistic
structure. However, in the case when the dimension of $H$ is not
less than 3, Uhlhorn was able to obtain the same conclusion as
Wigner.

In the last decades it has become quite clear that indefinite
inner product spaces are even more useful than definite ones in
describing several physical problems (see, for example, the
introduction in \cite{BMS}). This has raised the need to study
Wigner's theorem in the indefinite setting as well (see \cite{BMS}
and \cite{Broek}). Our paper \cite{MolCMP} was devoted to a
generalization of Wigner's original theorem for indefinite inner
product spaces. In the present paper we treat Uhlhorn's version in
that setting. Our approach here is different from that followed in
\cite{MolCMP}. Namely, it is based on a beautiful result of
Ovchinnikov \cite{Ovchinnikov} describing the automorphisms of the
poset of all idempotents on a separable Hilbert space of dimension
at least 3, which result can be regarded as a ``skew version'' of
the fundamental theorem of projective geometry. This result
enables us to use operator algebraic tools to attack the problem.
We note that this kind of machinery already proved effective in
our former works \cite{MolJMP,MolCMP2} where we obtained
Wigner-type results for different structures. We emphasize that in
the literature there does exist an Uhlhorn-type result on symmetry
transformations on indefinite inner product spaces. In fact, this
is due to Van den Broek \cite{Broek} (an application of his result
can be found in \cite{Broek2}, also see \cite{Broek3}). In that
paper he considered indefinite inner product spaces induced by
nonsingular self-adjoint operators on finite dimensional complex
Hilbert spaces. Moreover, in the proof of the main result he
basically followed the original idea of Uhlhorn. In the present
paper, we apply a completely different approach and obtain a much
more general result, namely, a result concerning indefinite inner
product spaces induced by any invertible bounded linear operator
on a real or complex Hilbert space of any dimension (not less than
3). Quantum logics on spaces with such a general indefinite metric
have been investigated by, for example, Matvejchuk in
\cite{Matvejchuk}. Our result will follow from the main theorem of
the paper, which describes the form of all bijective
transformations of the set of all rank-one idempotents on a Banach
space which preserve zero products in both directions.

If $X$ is a (real or complex) Banach space, then $B(X)$ stands for
the algebra of all bounded linear operators on $X$. An operator
$P\in B(X)$ is called an idempotent if $P^2=P$. The set of all
idempotents in $B(X)$ is denoted by $I(X)$ and $I_1(X)$ stands for
the set of all rank-one elements of $I(X)$.

Now, our main result reads as follows.

\begin{mainthm}\label{T:uhl2}
Let $X$ be a (real or complex) Banach space of dimension at least 3.
Let $\phi:I_1(X) \to I_1(X)$ be a bijective transformation with the
property that
\[
PQ=0 \Longleftrightarrow \phi(P)\phi(Q)=0
\]
for all $P,Q\in I_1(X)$.

If $X$ is real, then
there exists an invertible bounded linear operator $A:X\to X$ such
that $\phi$ is of the form
\begin{equation}\label{E:uhl1}
\phi(P)=APA^{-1} \qquad (P \in I_1(X)).
\end{equation}

If $X$ is complex and infinite dimensional, then there exists an
invertible bounded linear or conjugate-linear operator $A: X\to X$
such that $\phi$ is of the form \eqref{E:uhl1}.

If $X$ is complex and finite dimensional, then we can suppose that
our transformation $\phi$ acts on the space of $n\times n$ complex
matrices ($n=\dim X$). In this case there is a nonsigular matrix
$A\in M_n(\C)$ and a ring automorphism $h$ of $\C$ such that
$\phi$ is of the form
\begin{equation}\label{E:uhl3}
\phi(P)=A h(P) A^{-1} \qquad (P\in I_1(\C^n)).
\end{equation}
Here $h(P)$ denotes the matrix obtained from $P$ by applying $h$ to
every entry of it.
\end{mainthm}

Our main theorem can be summarized by saying that every bijective
transformation on $I_1(X)$ which preserves zero products in both
directions comes from a linear or conjugate-linear algebra
automorphism of $B(X)$ if $X$ is real or complex and infinite
dimensional, and it comes from a semilinear algebra automorphism
of $B(X)$ if $X$ is complex and finite dimensional. Replying to a
remark of the referee, we note that our result probably has no
serious physical meaning. This is because the poset of all
idempotents on a Banach space (the partial order among idempotents
is defined in Section 2.) does not form a lattice in general and
hence it is not a geometry or a logic in the sense of quantum
mechanics (see \cite{Varadarajan}). In fact, the poset of
idempotents is not to be confused with the lattice of subspaces of
a linear space as the idempotents are determined not by one but
two complementary subspaces. However, our main theorem will easily
imply our result Corollary~\ref{C:uhl2} generalizing Uhlhorn's
version of Wigner's theorem for indefinite inner product spaces
which statement we believe has serious physical meaning. On the
other hand, it will be clear from the proof presented that one can
readily get a very similar result as in our theorem for the form
of zero product preserving transformations on the set of rank-one
idempotents on different Banach spaces (also see the remark after
Corollary~\ref{C:uhl2}) which has an interesting mathematical
consequence. Namely, it implies that the real Banach spaces as
topological vector spaces are completely determined by the set of
their rank-one idempotents with the relation of zero product.

In our paper \cite{MolCMP} we presented a Wigner-type result for
pairs of ray transformations (\cite[Theorem 1]{MolCMP}) which
enabled us to generalize the result of Bracci, Morchio and
Strocchi in \cite{BMS} for indefinite inner product spaces
generated by any invertible bounded linear (not necessarily
self-adjoint) operator on a Hilbert space. Now, our main result
above can be applied to obtain the following corollary, which is a
Banach space analogue and hence a remarkable generalization (in
the mathematical sense) of the main result in \cite{MolCMP} that
was formulated for (complex) Hilbert spaces.

For the formulation of our corollary we need some concepts and
notation. Following the terminology of Uhlhorn, for any vector
$x\in X$, the set $\x$ of all nonzero scalar multiples of $x$ is
called the ray generated by $x$. The set of all rays in $X$ is
denoted by $\X$. The dual space of $X$ (that is the set of all
bounded linear functionals on $X$) is denoted by $X'$. For any
$x\in X, f\in X'$ we use the common and convenient notation $\la
x,f\ra$ for $f(x)$. We say that the rays $\x \in \X$ and $\f \in
\Xs$ are orthogonal to each other, in notation $\x \cdot \f=0$, if
we have $\la y, g\ra =0$ for all $y\in \x$ and $g\in \f$. The
Banach space adjoint of an operator $A\in B(X)$ is denoted by
$A'$. We extend the concept of adjoints also for conjugate-linear
operators. If $A$ is a bounded conjugate-linear operator on the
complex Banach space $X$, then its adjoint $A':X'\to X'$ (which is
also a bounded conjugate-linear operator) is defined by $A'
f=\overline{f\circ A}$ $(f\in X')$. If $X$ is a linear space over
$\K$ ($\K$ denotes the real or complex field) and $h$ is a ring
automorphism of $\K$, then the function $A:X \to X$ is called
$h$-semilinear if it is additive and $A(\lambda x)=h(\lambda) Ax$
holds for every $x\in X$ and $\lambda \in \K$. If $X$ is a finite
dimensional complex linear space and $h$ is a ring automorphism of
$\C$, then for any $h$-semilinear operator $A$, the adjoint $A'$
of $A$ is defined by $A'f=h^{-1}\circ f\circ A$ $(f\in X')$.
Clearly, $A':X'\to X'$ is an $h^{-1}$-semilinear operator.

After this preparation we can formulate our first corollary as
follows.

\begin{corollary}\label{C:uhl1}
Let $X$ be a (real or complex) Banach space of dimension not less than
3.
Let $T: \X \to \X$ and $S:\Xs \to \Xs$ be bijective transformations with
the property that
\[
T\x \cdot S\f =0 \quad \text{ if and only if } \quad
\x \cdot \f=0
\]
for every $\x \in \X$ and $\f \in \Xs$.

If $X$ is real, then
there exists an invertible bounded linear operator $A:X\to X$ such
that $T,S$ are of the forms
\begin{equation}\label{E:uhl6}
T\x ={\underline{Ax}} \quad \text{and} \quad
S\f={\underline{{A^{-1}}'f}} \qquad (0\neq x \in X, 0\neq f \in X').
\end{equation}

If $X$ is complex and infinite dimensional, then there exists an
invertible bounded linear or conjugate-linear operator $A: X\to X$
such that $T,S$ are of the forms \eqref{E:uhl6}.

If $X$ is complex and finite dimensional, then
there exist a ring automorphism $h$ of $\C$ and
an invertible $h$-semilinear operator $A: X\to X$ such that $T,S$ are
of the forms \eqref{E:uhl6}.

The operator $A$ above is unique up to multiplication by a scalar.
\end{corollary}

Finally, as a consequence of Corollary~\ref{C:uhl1}, we shall
present our Uhlhorn-type version of Wigner's theorem for
indefinite inner product spaces that was promised in the abstract.
As mentioned above, our result is a far-reaching generalization of
the main result in \cite{Broek}, where a similar assertion in the
particular case when $H$ is finite dimensional and the generating
invertible operator $\eta$ is self-adjoint was presented.

Let $\eta$ be an invertible bounded linear operator on a Hilbert space
$H$. Denote by $(x, y)_\eta$ the quantity $\la \eta x,y\ra$ $(x, y\in
H)$. We write $\x \cdot_\eta \y=0$ if $\la \eta x_0, y_0\ra=0$ holds for
every $x_0 \in \x$ and $y_0\in \y$.
The ray transformation $T:\H \to \H$ is called a symmetry transformation
on the indefinite inner product space $\H$ generated by $\eta$ if
\[
T\x \cdot_\eta T\y=0 \quad \Longleftrightarrow \quad \x \cdot_\eta \y=0
\]
for all $\x, \y\in \H$. We say that the transformation $T:\H \to
\H$ is induced by the invertible linear or conjugate-linear
operator $U: H\to H$ if $T\x=\underline{Ux}$ for every $0\neq x\in
H$.

\begin{corollary}\label{C:uhl2}
Let $H$ be a (real or complex) Hilbert space of dimension not less than
3 and let $\eta \in B(H)$ be invertible. Suppose that $T:\H \to \H$ is a
bijective transformation with the property that
\[
T\x \cdot_\eta T\y =0 \quad \text{ if and only if } \quad
\x \cdot_\eta \y=0
\]
holds for every $\x,\y \in \H$.

If $H$ is real, then $T$ is induced by an invertible bounded
linear operator $U$ on $H$. Similarly, if $H$ is complex, then $T$
is induced by an invertible bounded linear or conjugate-linear
operator $U$ on $H$.

The operator $U$ inducing $T$ is unique up to multiplication by a
scalar.

If $H$ is real, then the invertible bounded linear operator
$U:H\to H$ induces a symmety transformation on $\H$ if and only if
\[
(Ux,Uy)_\eta =c (x,y)_\eta \qquad (x,y \in H)
\]
holds for some constant $c\in \R$.

If $H$ is complex, then the invertible bounded linear operator
$U:H\to H$ induces a symmetry transformation on $\H$ if and only if
\[
(Ux,Uy)_\eta =c (x,y)_\eta \qquad (x,y \in H)
\]
holds for some constant $c\in \C$. Similarly,
the invertible bounded conjugate-linear operator
$U:H\to H$ induces a symmetry transformation on $\H$ if and only if
\[
(Ux,Uy)_\eta =d (y,x)_{\eta^*}  \qquad (x,y \in H)
\]
holds for some constant $d\in \C$.
Here, $\eta^*$ denotes the Hilbert space adjoint of $\eta$.
\end{corollary}

\begin{remark}
Observe that in contrast with the Main Theorem and
Corollary~\ref{C:uhl1}, in Corollary~\ref{C:uhl2} above general
semilinear operators do not appear.

In Uhlhorn's paper \cite{Uhlhorn} it was mentioned that, for
physical reasons, one should consider ray transformations between
different spaces. It will be clear from the proofs below that one
can generalize our result in that direction easily.

We should point out that, as will be clear from their proofs, in
Corollary~\ref{C:uhl1} and Corollary~\ref{C:uhl2} there is in fact
no need to assume the injectivity of the transformations $T,S$. We
have posed this condition only for the sake of ``symmetricity''.

Finally, we note that we are convinced that our result could
somehow be extended for the case of quaternionic Hilbert spaces,
which have also been proved to be important in the applications of
mathematics in certain physical problems. The first step in this
direction could be an extension of Ovchinnikov's result for that
case. However, we leave the whole (we believe challenging) problem
open.
\end{remark}

\section{Proofs}

In the proofs we need some additional notation and definitions.

Let $X$ be a (real or complex) Banach space.
The ideal of all finite rank operators in
$B(X)$ is denoted by $F(X)$.
Two idempotents $P,Q$ in $B(X)$ are said to be
(algebraically)
orthogonal if $PQ=QP=0$. There is a natural partial order on
$I(X)$. Namely, for any $P,Q\in I(X)$
we write $P\leq Q$ if $PQ=QP=P$. Clearly, $P\leq Q$ holds if and
only if the range $\rng P$ of $P$ is a subset of the range of $Q$ and
the kernel $\ker P$ of $P$ contains the kernel of $Q$.
The symbol $I_f(X)$ stands for the collection
of all finite rank idempotents in $B(X)$.
The natural embedding of $X$ into its second dual $X''$ is denoted by
$\kappa$. If $x\in X$ and $f\in X'$, then $x\ot f$ stands for the
operator (of rank at most 1) defined by
\[
(x\ot f)(z)=\la z, f\ra x \qquad (z\in X).
\]
Clearly, $x\ot f$ is a rank-one idempotent if and only if $\la
x,f\ra=1$. It is easy to see that the elements of $F(X)$ are
exactly the operators $A\in B(X)$ which can be written as finite
sums of the form
\begin{equation}\label{E:uhl11}
A=\sum_i x_i \ot f_i
\end{equation}
with $x_1, \ldots, x_n\in X$ and $f_1, \ldots, f_n\in X'$. Using
this representation, the trace of $A$ is defined by
\[
\tr A=\sum_i \la x_i, f_i\ra.
\]
It is known that $\tr A$ is well defined, that is, it does not
depend on the particular representation \eqref{E:uhl11} of $A$.
Denote by $M_n(\K)$ the algebra of all $n\times n$ matrices with
entries in $\K$.

In the proof of our main result we shall need the following lemma.

\begin{lemma}
For any $P_1, P_2\in I_f(X)$ there
exists a $P\in I_f(X)$ such that
$P_1,P_2 \leq P$.
\end{lemma}

\begin{proof}
The assertion will follow from the following observation. Let $M,N
\subset X$ be closed subspaces. Suppose that $M$ is of finite
codimension and $N$ is of finite dimension. Then there exists an
idempotent $P\in I_f(X)$ such that $\ker P\subset M$ and $\rng
P\supset N$. Indeed, since every finite-dimensional subspace of a
Banach space is complemented, we can find a closed subspace $K$ in
$X$ such that $K\oplus (M\cap N)=M$. Since the sum of a closed and
a finite dimensional subspace is closed, it follows that $M+N$ is
closed and has finite codimension. So, there is a
finite-dimensional subspace $L$ in $X$ such that $(M+N)\oplus
L=X$. We clearly have
\[
K\oplus(N\oplus L)=X.
\]
Now, there exists an idempotent $P\in I_f(X)$ such that
$\ker P=K$ and $\rng P=N\oplus L$. This verifies our observation.

If $P_1, P_2\in I_f(X)$, then
$\ker P_1\cap \ker P_2$ is of finite corank and $\rng P_1 +\rng
P_2$ is of finite rank. Now, the
idempotent $P\in I_f(X)$ obtained
according to the observation above
clearly has the property that $P_1,P_2\leq P$.
This completes the proof.
\end{proof}

\begin{proof}[Proof of the Main Theorem]
We first extend $\phi$ to the set $I_f(X)$ of all finite rank
idempotents in $B(X)$. If $0\neq P\in I_f(X)$, then there are
mutually (algebraically) orthogonal rank-one idempotents $P_1,
\ldots, P_n\in B(X)$ such that $P=\sum_i P_i$. Clearly,
$\phi(P_1), \ldots, \phi(P_n)$ are also mutually orthogonal
rank-one idempotents. Let us define
\[
\tphi(P)=\sum_i \phi(P_i).
\]
We have to show that $\tphi$ is well defined. In order to do this,
let $Q_1, \ldots, Q_n\in B(X)$ be mutually orthogonal rank-one
idempotents with sum $P$. Pick any $R\in I_1(X)$. We have
\begin{equation*}
\begin{gathered}
(\sum_i \phi(P_i))\phi(R)=0 \Longleftrightarrow
\phi(P_i)\phi(R)=0 \enskip (i=1, \ldots, n) \Longleftrightarrow\\
P_iR=0 \enskip (i=1, \ldots, n) \Longleftrightarrow
(\sum_iP_i)R=0.
\end{gathered}
\end{equation*}
Similarly, we obtain
\[
(\sum_i \phi(Q_i))\phi(R)=0 \Longleftrightarrow
(\sum_iQ_i)R=0.
\]
Since $\sum_i P_i =\sum_i Q_i$,
these imply that
\[
(\sum_i \phi(P_i))\phi(R)=0 \Longleftrightarrow
(\sum_i \phi(Q_i))\phi(R)=0.
\]
As $\phi(R)$ runs through the set $I_1(X)$, we deduce that the
kernels of the idempotents $\sum_i \phi(P_i)$ and $\sum_i
\phi(Q_i)$ are the same. A similar argument shows that the ranges
of these two idempotents are also equal. Therefore, we have
\[
\sum_i \phi(P_i)=\sum_i \phi(Q_i).
\]
This shows that the transformation $\tphi$ is well defined. It is
now easy to verify that $\tphi :I_f(X) \to I_f(X)$ is a bijection
which preserves the order, the orthogonality and the rank in both
directions. In fact, only the injectivity is not trivial but it
follows from an argument quite similar to the one proving $\tphi$
is well defined.

Pick a finite rank idempotent $P_0 \in B(X)$ whose rank is at least 3.
Consider the set $I_{P_0}(X)$ of all idempotents $P\in B(X)$ for which
$P\leq P_0$. Let $M=\ker P_0$ and $N=\rng P_0$. We have $M\oplus N=X$.
Denote by $B(X,M,N)$ the set of all operators $A$ in $B(X)$ for which
$A(N)\subset N$ and $A(M)=\{ 0\}$.
Clearly, we have $I_{P_0}(X) \subset B(X,M,N)$.
Considering the transformation
$A\longmapsto A_{|N}$ we get an algebra isomorphism from $B(X,M,N)$ onto
$B(N)$. Moreover, $B(N)$ is obviously isomorphic to $M_n(\K)$.
Denote the so-obtained algebra isomorphism from $B(X,M,N)$ onto
$M_n(\K)$ by $\psi$.
Similarly, we have an algebra isomorphism $\psi'$ from
$B(X, \ker \phi(P_0), \rng \phi(P_0))$ onto $M_n(\K)$.
Therefore, the transformation $P\mapsto \Psi(P)=\psi'(\tphi
(\psi^{-1}(P)))$ is a bijection of the set of all idempotents in
$M_n(\K)$ which preserves the order $\leq$ in both directions.
The form of all such transformations is described on p. 186 in
\cite{Ovchinnikov}. In particular,
it follows from that form that there is a ring-automorphism
$h_{P_0}$ of $\K$ such that
\[
\tr \Psi(P)\Psi(Q)=h_{P_0}(\tr PQ)
\]
holds for all idempotents $P,Q$ in $M_n(\K)$. Since $\psi, \psi'$
are algebra isomorphisms, it follows that they preserve rank-one
idempotents. This implies that $\psi, \psi'$ preserve the traces
of rank-one operators, from which we conclude that they are
generally trace-preserving. It follows that
\begin{equation}\label{E:uhl100}
\tr \tphi(P)\tphi(Q)=h_{P_0}(\tr PQ) \qquad (P,Q\in I_{P_0}(X)).
\end{equation}
We claim that in fact $h_{P_0}$ does not depend on $P_0$. Indeed, let
$P_1\in I_f(X)$ be such that $P_0 \leq
P_1$. Considering the corresponding ring automorphism $h_{P_1}$ of
$\K$, by \eqref{E:uhl100} we get that
\[
h_{P_0}( \tr PQ)=h_{P_1}( \tr PQ)
\]
holds for every $P,Q \in I_{P_0}(X)$. Clearly, $\tr PQ$ runs through
$\K$
as $P,Q$ run through $I_{P_0}(X)$. This shows that $h_{P_0}=h_{P_1}$.
Since for any two finite rank idempotents there is a finite rank
idempotent majorizing both of them (this is just
the content of our Lemma), we have the independence of $h_{P_0}$ from
$P_0$. Therefore, there exists a ring automorphism $h$ of $\K$ such that
\begin{equation}\label{E:uhl13}
\tr \tphi (P)\tphi(Q)=h(\tr PQ) \qquad (P,Q\in I_f(X)).
\end{equation}
We now extend $\tphi$ from $I_f(X)$ onto $F(X)$. For any $P_1, \ldots,
P_n \in I_f(X)$ and $\lambda_1, \ldots, \lambda_n\in \K$ we define
\[
\Phi(\sum_i \lambda_i P_i)=\sum_i h(\lambda_i) \tphi(P_i).
\]
We have to show that $\Phi$ is well defined. Let $Q_1, \ldots,
Q_m\in I_f(X)$ and $\mu_1, \ldots, \mu_m\in \K$ be such that
\[
\sum_i \lambda_i P_i=\sum_j \mu_j Q_j.
\]
It follows that
\[
\sum_i \lambda_i P_iR=\sum_j \mu_j Q_jR
\]
holds for every $R\in I_f(X)$. Taking traces we obtain
\[
\sum_i \lambda_i \tr P_iR=\sum_j \mu_j \tr Q_jR.
\]
By \eqref{E:uhl13} it follows that
\[
\sum_i \lambda_i h^{-1}(\tr \tphi(P_i)\tphi(R))=
\sum_j \mu_j h^{-1}(\tr \tphi(Q_j)\tphi(R)).
\]
This implies that
\[
h^{-1}(\sum_i h(\lambda_i) \tr \tphi(P_i)\tphi(R))=
h^{-1}(\sum_j h(\mu_j) \tr \tphi(Q_j)\tphi(R)),
\]
that is,
\[
h^{-1}(\tr (\sum_i h(\lambda_i) \tphi(P_i))\tphi(R)))=
h^{-1}(\tr (\sum_j h(\mu_j)  \tphi(Q_j))\tphi(R))).
\]
This gives
\[
\tr (\sum_i h(\lambda_i) \tphi(P_i))\tphi(R)=
\tr (\sum_j h(\mu_j)  \tphi(Q_j))\tphi(R).
\]
Since $\tphi(R)$ runs through the set $I_f(X)$, we obtain
\[
\sum_i h(\lambda_i) \tphi(P_i)=
\sum_j h(\mu_j)  \tphi(Q_j).
\]
Therefore, $\Phi$ is well defined.
Since the finite rank idempotents linearly generate $F(X)$,
it follows that $\Phi$ is a surjective
$h$-semilinear transformation on $F(X)$
which preserves the rank-one idempotents and their
linear spans. We can now apply a result of Omladi\v c and \v Semrl
describing the form of all such transformations.
In fact, if, for example, $X$ is real, then by
\cite[Main Result]{Semrl} either there exists an invertible bounded
linear operator $A:X\to X$ such that
\begin{equation}\label{E:uhl51}
\phi(P)=APA^{-1} \qquad (P\in I_1(X))
\end{equation}
or there exists an invertible bounded linear operator
$B:X'\to X$ such that
\begin{equation*}\label{E:uhl52}
\phi(P)=BP'B^{-1} \qquad (P\in I_1(X)).
\end{equation*}
If we had this second possibility, then we would get that
\[
\phi(P)\phi(Q)=0 \Longleftrightarrow
BP'Q'B^{-1}=0    \Longleftrightarrow
P'Q'=0           \Longleftrightarrow
QP=0
\]
for every $P,Q\in I_1(X)$. On the other hand, we know that
\[
\phi(P)\phi(Q)=0 \Longleftrightarrow
PQ=0.
\]
So, we would have
\[
PQ=0
\Longleftrightarrow
QP=0
\]
for every $P,Q\in I_1(X)$, which is an obvious contradiction.
Therefore, $\phi$ is of the form \eqref{E:uhl51}.

If $X$ is complex, then one can argue in a very similar way referring to
\cite[Main Result]{Semrl} again (in the infinite dimensional case) or
to \cite[Theorem 4.5]{Semrl} (in the finite dimensional case).
The proof is complete.
\end{proof}

\begin{proof}[Proof of Corollary~\ref{C:uhl1}]
We define a bijective transformation $\phi:I_1(X) \to I_1(X)$
which preserves zero products in both directions.

First, for every $0\neq x\in X$ pick a vector from the ray $T\x$.
In that way we get a transformation, which
will be denoted by the same symbol $T$, from $X\setminus \{ 0\}$ into
itself with the
property that for every vector $0\neq y\in X$, there exists a vector
$0\neq x\in X$ such
that $y=\lambda Tx$ for some nonzero scalar $\lambda\in \K$. We do the
same with the other transformation $S$. Clearly, we have
\begin{equation}\label{E:uhl10}
\la Tx, Sf\ra =0 \quad
\text{ if and only if } \quad
\la x, f\ra =0
\end{equation}
for every nonzero $x\in X$ and nonzero $f\in X'$.

Let $x\in X$ and $f\in X'$ be such that $\la x,f\ra \neq 0$.
Define
\[
\phi\biggl(\frac{1}{\la x, f\ra} x \ot f\biggr)=
\frac{1}{\la Tx, Sf\ra} Tx \ot Sf.
\]
We show that $\phi$ is well defined. Let $x_0\in X$ and $f_0\in
X'$ be such that $\la x_0,f_0\ra \neq 0$ and suppose that
\[
\frac{1}{\la x, f\ra} x \ot f=
\frac{1}{\la x_0, f_0\ra} x_0 \ot f_0.
\]
This implies that $x,x_0$ belong to the same ray in $X$ and the
same holds true for $f,f_0$ in $X'$. Consequently, $Tx, Tx_0$ and
$Sf,Sf_0$ generate equal rays in $X$ and $X'$, respectively.
Therefore, the ranges and the kernels of the idempotents
$\frac{1}{\la Tx, Sf\ra} Tx \ot Sf$ and $\frac{1}{\la
Tx_0,Sf_0\ra} Tx_0 \ot Sf_0$ are equal, which implies the equality
of these two idempotents. Hence, we obtain that $\phi$ is well
defined.

By the "almost surjectivity" property of the vector-vector
transformations $T,S$ we obtain the surjectivity of $\phi$. The
injectivity of $\phi$ can be proved by an argument like the one we
used to prove $\phi$ is well defined. The transformation $\phi$
preserves zero products in both directions, which is a consequence
of \eqref{E:uhl10}.

Now, we can apply our main theorem.
Suppose first that $X$ is real. Then our transformation $\phi$ is of
the
form \eqref{E:uhl1} with some invertible bounded linear operator $A$ on
$X$. If $x\in X$ and $f\in X'$ are such that $\la x,f\ra\neq 0$, then
from the equality
\begin{equation}\label{E:uhl120}
\begin{gathered}
\frac{1}{\la Tx,Sf\ra}Tx\ot Sf=
\phi\biggl(\frac{1}{\la x,f\ra}x\ot f\biggr)=\\
A \cdot \frac{1}{\la x,f\ra}x\ot f\cdot A^{-1}=
\biggl(\frac{1}{\la x,f\ra}Ax\biggr)\ot ({A^{-1}}'f)
\end{gathered}
\end{equation}
we deduce that $Tx$ is a scalar multiple of $Ax$ and
$Sf$ is a scalar multiple of ${A^{-1}}'f$.
This gives us that $T\x=\underline{Ax}$ and $S\f=\underline{{A^{-1}}'
f}$.

If $X$ is complex infinite dimensional, then one can argue in a very
similar way.

Finally, let $X$ be complex and finite dimensional.
In that case there exist a ring automorphism $h$ of $\C$ and an
invertible $h$-semilinear operator $A:X\to X$ such that
$\phi$ is of the form
\begin{equation*}
\phi(P)=A P A^{-1} \qquad (P\in I_1(X)).
\end{equation*}
This comes from a rewriting of the form
\eqref{E:uhl3} appearing in the formulation of our main theorem.
Now, one can easily verify that we have the following equality very
similar to \eqref{E:uhl120}:
\begin{equation*}
\frac{1}{\la Tx,Sf\ra}Tx\ot Sf=
\biggl(\frac{1}{h(\la x,f\ra)}Ax\biggr)\ot ({A^{-1}}'f).
\end{equation*}
This yields
$T\x=\underline{Ax}$ and $S\f=\underline{{A^{-1}}'
f}$ $(\x \in \X, \f \in \X')$.

The assertion concerning essential uniqueness is a consequence of the
following easy fact whose proof requires only elementary linear algebra.
If $A,B$ are semilinear operators on a vector space $Y$ over $\K$ with
ranks at least 2 such that $Ay,By$ are linearly dependent for every
$y\in Y$, then $A,B$ are linearly dependent.
This completes the proof of Corollary~\ref{C:uhl1}.
\end{proof}

\begin{proof}[Proof of Corollary~\ref{C:uhl2}]
Just as in the proof of Corollary~\ref{C:uhl1}, we can define an
"almost surjective" transformation (that is, one that has values
in every ray) on the underlying Hilbert space $H$, denoted by the
same symbol $T$, such that
\[
\la \eta Tx, Ty\ra =0 \quad \text{if and only if} \quad
\la \eta x, y\ra =0 \quad (x, y\in H\setminus \{ 0\}).
\]
We can rewrite this equivalence first as
\[
\la \eta T \eta^{-1} x, Ty\ra =0 \quad \text{if and only if} \quad
\la x, y\ra =0 \quad (x, y\in H\setminus \{ 0\})
\]
and next as
\[
\la Tx, \eta T \eta^{-1} y\ra =0 \quad \text{if and only if} \quad
\la x, y\ra =0 \quad (x, y\in H\setminus \{ 0\}).
\]
Now, we apply Corollary~\ref{C:uhl1}. To be honest, we should
point out that although that result is formulated for Banach
spaces and hence dual spaces and Banach space adjoints of
operators appear there, the very same argument can be applied to
conclude that our present transformation $T$ is generated by some
invertible operator $U$ on $H$. We learn from
Corollary~\ref{C:uhl1} that $U$ is linear if $H$ is real, it is
either linear or conjugate-linear if $H$ is complex infinite
dimensional and, finally, $U$ is semilinear if $H$ is complex
finite dimensional. From the proof of the remaining part of our
corollary it will be clear that this general semilinear case in
fact does not occur.

The essential uniqueness of $U$ can be verified as in the proof of
Corollary~\ref{C:uhl1}. As for the third part of the statement, we
present the proof only in the complex finite-dimensional case. In
all other cases one can argue in a quite similar way. So, let $h$
be a ring automorphism of $\C$. Suppose that the invertible
$h$-semilinear operator $U: H\to H$ induces a symmetry
transformation. Then we have
\[
\la \eta Ux, Uy\ra =0 \quad \Longleftrightarrow \quad
\la \eta x,y\ra =0
\]
for every $x,y \in H$.
This implies that
\[
h^{-1}(\la \eta Ux, Uy\ra) =0
\quad \Longleftrightarrow \quad
\la \eta x,y\ra =0 \qquad (x,y \in H).
\]
If we fix $y\in H$, then the functions $x\mapsto h^{-1}(\la \eta
Ux, Uy\ra)$ and $x\mapsto \la \eta x,y\ra$ are linear functionals
with the same kernel. We deduce that these functionals differ only
by a scalar multiple. Hence, there exists a $c(y)\in \C$ such that
\begin{equation}\label{E:uhl20}
h^{-1}(\la \eta Ux, Uy\ra) =
c(y) \la \eta x,y\ra
\end{equation}
for every $x,y \in H$. Similarly, for every $x\in H$ there exists a
scalar $d(x)\in \C$ such that
\[
h^{-1}(\la Uy, \eta Ux\ra)=
d(x)\la y,\eta x\ra \qquad (x,y\in H).
\]
Defining $g:\C \to \C$ by
$g(\lambda)=\overline{h(\overline{\lambda})}$ $(\lambda\in \C)$,
we can write this last equality as
\begin{equation}\label{E:uhl21}
g^{-1}(\la \eta Ux, Uy\ra)=\overline{d(x)}\la \eta x,y\ra
\qquad (x,y\in H).
\end{equation}
It follows from \eqref{E:uhl20} and \eqref{E:uhl21} that
\[
\la \eta Ux, Uy\ra=C(y) h(\la \eta x,y\ra)
\enskip \text{ and } \enskip
\la \eta Ux, Uy\ra=D(x) g(\la \eta x,y\ra)
\]
for every $x,y\in H$, where $C,D$ are complex-valued functions on
$H$. We then have
\[
C(y) h(\la \eta x,y\ra)=
D(x) g(\la \eta x,y\ra)
\]
for every $x,y\in H$. It is easy to see that $C,D$ are in fact
constant functions. Indeed, pick any $y_1,y_2\in H$ which are
linearly independent. Then we have $x,z\in H$ such that $y_1=\eta
x$, $z\perp \eta x$ and $y_2=\eta x+z$. Since $\la \eta x,y_1\ra
=\la \eta x,y_2\ra$, it follows from the equality above that
$C(y_1)=C(y_2)$. In case $y_1, y_2\in H\setminus \{ 0\}$ are
linearly dependent, we can choose $y_3\in H$ such that $y_1,y_3$
and $y_2,y_3$ are both linearly independent and we get
$C(y_1)=C(y_2)$. Since $C(0)$ does not count, we obtain that $C$
is really constant. A similar argument applies to $D$. It follows
that we have constants $C,D\in \C$ such that
\[
\la \eta Ux, Uy\ra =
C h(\la \eta x,y\ra)
\]
and
\[
\la \eta Ux, Uy\ra=D \overline{h(\overline{\la \eta x,y\ra})}.
\]
Since these hold for every $x,y\in H$ and we have $h(1)=1$, it
follows that $C=D$. This implies that $h$ is self-adjoint in the
sense that $h(\overline{\lambda})=\overline{h(\lambda)}$ $(\lambda
\in \C)$. It is well known that the only ring automorphisms of
$\C$ with this property are the identity and the conjugation. In
fact, this is an easy consequence of the fact that the only ring
automorphism of $\R$ is the identity. It now follows that either
$U$ is linear and we have
\begin{equation}\label{E:uhl23}
(Ux,Uy)_\eta =C (x,y)_\eta \qquad (x,y \in H)
\end{equation}
or $U$ is conjugate-linear and we have
\begin{equation}\label{E:uhl24}
(Ux,Uy)_\eta =C (y,x)_{\eta^*}  \qquad (x,y \in H).
\end{equation}
It is obvious that if $U:H\to H$ is either an invertible linear
operator on $H$ such that \eqref{E:uhl23} holds or an invertible
conjugate-linear operator such that \eqref{E:uhl24} holds, then
$U$ induces a symmetry transformation.

The remaining part of the proof can be carried out in a similar,
but simpler, way.
\end{proof}

\section*{Acknowledgements}

The author is grateful to the referee whose kind remarks helped
him to make the presentation of the paper clearer.

\bibliographystyle{amsplain}

\end{document}